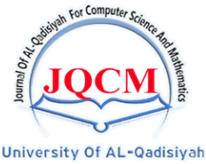

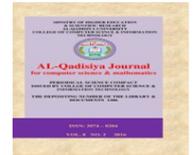

# On a Subclass of Meromorphic Univalent Functions Involving Hypergeometric Function

*Mazin Sh.Mahmoud[a], Abdul Rahman S.Juma[b], Raheam A. Mansor Al-Saphory[c]*

[a] Department of Mathematics, Tikrit University Tikrit, Iraq. Email: maz2004in@yahoo.com

[b] Department of Mathematics, University of Anbar, Ramadi, Iraq. Email: : dr-juma@hotmail.com

[c] Department of Mathematics, Tikrit University, Tikrit, Iraq, Email: saphory@hotmail.com

A B S T R A C T

The main object of the present paper is to, introduce the. class of meromorphic univalent functions Involving! hypergeomatrc function .We obtain~ some interesting geometric properties according to coefficient inequality, growth and distortion theorems, radii of starlikeness and convexity for the" functions in our subclass.

MSC.

Corresponding author Mazin Sh.Mahmoud

Email addresses: maz2004in@yahoo.com

Communicated by Qusuay Hatim Egaar

## 1 . Introduction

Let $\sum$ denoted be the class of functions of the form

$$f(z) = \frac{1}{z} + \sum_{k=1}^{\infty} a_k z^k \qquad (1)$$

which are analytic and meromorphic univalent in punctured unit disk $Ụ^* = \{z \in \mathbb{C} : 0 < |z| < 1\}$.

-A function $f \in \sum$ is meromorphic starlike of order $\alpha$, $(0 \leq \alpha < 1)$ if $R\left(\frac{zf'(z)}{f(z)}\right) > \alpha$, $(z \in Ụ^*)$.

The class of all such function is denoted by $\sum^*(\alpha)$. A function $f \in \sum$ is meromorphic convex of order $\alpha$, $(0 \leq \alpha < 1)$
- R $\left(1+\frac{zf''(z)}{f'(z)}\right) > \alpha$, $(z \in Ụ^*)$. Let $\sum_q$ be the class of function $f \in \sum$, if

with $a_k \geq 0$. The subclass of $\sum_q$ consisting of starlike functions of order α is denoted by $\sum_q^*(\alpha)$ , and convex functions of order $\alpha$ by $\sum_q^k(\alpha)$ . Various subclasses of $\sum$ have been defined and studied by varions authors see [1 , 2 , 3 , 4 , 5, 6 , 7 ,10 , 11,12].

For function $f(z)$ given by (1) and $g(z) = \frac{1}{z} + \sum_{k=1}^{\infty} b_k z^k$ ,we define theHadamard product or( convolution) $f$ and $g$ by

$$f * g = \frac{1}{z} + \sum_{k=1}^{\infty} a_k b_k z^k .$$

For positive real parameters ($\alpha_1, A_1, \ldots, \alpha_\ell, A_\ell, \beta_1, B_1, \ldots, \beta_p, B_p$)
($\ell, p \in \mathbb{N} = \{1, 2, \ldots\}$) such that

$1 + \sum_{n=1}^{p} B_n - \sum_{n=1}^{\ell} A_n \geq 0 , (z \in \mathbb{U}^*)$     The Wright generalized hypergeometric function

$$= {}_\ell\Psi_p[(\alpha_t, A_t)_{1,\ell}, (\beta_t, B_t)_{1,p}; Z], \Psi_p[(\alpha_1, A_1), \ldots, (\alpha_\ell, A_\ell); (\beta_1, B_1), \ldots, (\beta_p, B_p); Z]_\ell$$

assigned by

${}_\ell\Psi_p[(\alpha_t, A_t)_{1,\ell}, (\beta_t, B_t)_{1,p}; Z] = \sum_{n=0}^{\infty} \{\prod_{t=0}^{\ell} \Gamma(\alpha_1 + nA_t)\} \{\prod_{t=0}^{p} \Gamma(\beta_t + n B_t)\}^{-1} \frac{z^n}{n!}$ .

If $A_t = 1, (t = 1, 2, 3, \ldots \ell)$ and $B_t = 1, (t = 1, 2, 3, \ldots, P)$, then

$\Omega_\ell \Psi_p[(\alpha_t A_t)_{1,\ell}, (\beta_t, B_t)_{1,p}; z] \equiv F_p(\alpha_1, \ldots, \alpha_\ell, \beta_1, \ldots, \beta_m; z)$

$$= \sum_{n=0}^{\infty} \frac{(\alpha_1)_n \ldots (\alpha_1)_n z^n}{(\beta_1)_n \ldots (\beta_p)_n n!}$$

($\ell \leq p + 1; \ell, p \in \mathbb{N}_0 = \mathbb{N} = \{0, 1, 2, 3, \ldots\}; Z \in \mathbb{U}$).

That is the generalized hypergeometric function (see[8]). Here $(\alpha_k)$ is the Pochammer symbol and $\Omega = (\prod_{t=0}^{\ell} \Gamma(\alpha_t))^{-1} (\prod_{t=0}^{p} \Gamma(\beta_t))$.

When assign the generalized hypergeometric function , we take a Linear operator

$W[(\alpha_1, A_t)_{1,\ell}, (\beta_t, B_t)_{1,P}]: \sum_q \quad \sum_q . \longrightarrow$

$W[(\alpha_t, A_t)_{1,\ell}, (\beta_t, B_t)_{1,P}] f(z) = z^{-1} \Omega_l \Psi_P[(\alpha_t, A_t)_{1,\ell} (\beta_t, B_t)_{1,P}; z] * f(z)$     (2)

for convenience ,we denote $W[(\alpha_t, A_t)_{1,\ell}, (\beta_t, B_t)_{1,p}]$ by $W[\alpha_1]$.
If $f$ has the from (1) then we obtain

$$W[\alpha_1] f(z) = \frac{1}{z} + \sum_{k=1}^{\infty} \sigma_k(\alpha_1) a_k z^k , \qquad (3)$$

where

$$\sigma_k(\alpha_1) = \frac{\Omega \Gamma(\alpha_1 + A_1(k+1)) \ldots \Gamma(\alpha_\ell + A_\ell(k+1))}{(k+1)! \, r(\beta_1 + B_1(k+1)) \ldots \Gamma(\beta_\ell + B_\ell(k+1))} . \qquad (4)$$

**Definition1.1**: A subclass of $\Sigma_q$ by utilizing operator $W[\alpha_1]$ we let $V(\alpha,\eta)$ denote a subclass of $\Sigma_q$ consisting of function in (1) satisfying the condition

$$\left| \frac{\frac{z(W(\alpha_1)f(z))''}{(W(\alpha_1)f(z))'} + 2}{\frac{z(W(\alpha_1)f(z))''}{(W(\alpha_1)f(z))'} + 2\alpha} \right| < \eta. \qquad (5)$$

$0 < \alpha < 1$, $0 < \eta \leq 1$ and $A_t = 1$ (t=1,2,3,…), $B_t = 1$ (t=1,2,3,…) Where
Now we must prove the Coefficient Inequality

## 2. Coeffcient Inequality

**Theorem 2.1**: $f$ is a function defined by (1) in the class $V(\alpha,\eta)$, if and only if

$$\sum_{k=1}^{\infty} |\sigma_k(\alpha_1)|[k(1+\eta) + (1+\eta(2\alpha-1))]a_k \leq 2\eta(1-\alpha). \qquad (6)$$

*The result is sharp"*

**Proof**: Let the inequity (6) holds true and let $|z|=1$ by (5). Then we get

$$\left| \frac{z(W(\alpha_1)f(z))''}{(W(\alpha_1)f(z))'} + 2 \right| - \eta \left| \frac{z(W(\alpha_1)f(z))''}{(W(\alpha_1)f(z))'} + 2\alpha \right| < 0,$$

$|z(W(\alpha_1)f(z))'' + 2(W(\alpha_1)f(z))'| - \eta|z(W(\alpha_1)f(z))'' + 2\alpha(W(\alpha_1)f(z))'|,$

and by utilizing (3) we have

$$(W(\alpha_1)f(z))' = \frac{-1}{z^2} + \sum_{k=1}^{\infty} n\,\sigma_k(\alpha_1)a_k z^{k-1},$$

$$(W(\alpha_1)f(z))'' = \frac{+2}{z^3} + \sum_{k=1}^{\infty} k(k-1)\,\sigma_k(\alpha_1)a_k z^{k-2},$$

$$z(w(\alpha_1)f(z))'' = \frac{2}{z^2} + \sum_{k=1}^{\infty} k(k-1)\sigma_k(\alpha_1)a_k z^{k-1},$$

$= \left| \frac{2}{z^2} + \sum_{k=1}^{\infty} k(k-1)\sigma_k(\alpha_1)a_k z^{k-1} - \frac{2}{z^2} + 2\sum_{k=1}^{\infty} k\sigma_k(\alpha_1)a_k z^{k-1} \right|$
$-\eta \left| \frac{2}{z^2} + \sum_{k=1}^{\infty} k(k-1)\sigma_k(\alpha_1)a_k z^{k-1} - \frac{2\alpha}{z^2} + 2\alpha \sum_{k=1}^{\infty} k\sigma_k(\alpha_1)a_k z^{k-1} \right|$

$= \left| \sum_{k=1}^{\infty} k\sigma_k(\alpha_1)a_k z^{k-1}(k-1+2) \right| - \eta \left| \frac{2}{z^2}(1-\alpha) + \sum_{k=1}^{\infty} k\sigma_k(\alpha_1)a_k z^{k-1})(k-1+2\alpha) \right|$

$$\leq \sum_{k=1}^{\infty} k\,\sigma_k(\alpha_1)a_k |z|^{k-1}(k+1) \frac{-2\eta}{|z|^2}(1-\alpha) + \eta \sum_{k=1}^{\infty} k\,\sigma_k(\alpha_1)a_k |z|^{k-1}(k-1+2\alpha)$$

$$\leq \sum_{k=1}^{\infty} |\sigma_k(\alpha_1)|[k(1+\eta) + (1+\eta(2\alpha-1)]a_k - 2\eta(1-\alpha) \leq 0.$$

Therefore, by the maximum modules theorem we have $f \in V(\alpha,\eta)$,
Conversely, suppose $f \in V(\alpha,\eta)$, then

$$\left| \frac{\frac{z(W(\alpha_1)f(z))''}{W(\alpha_1)f(z))'} + 2}{\frac{z(W(\alpha_1)f(z))''}{(W(\alpha_1)f(z))'} + 2\alpha} \right| < \eta,$$

$$\left| \frac{\sum_{k=1}^{\infty}(k+1)|\sigma_k(\alpha_1)|a_k z^{k-1}}{\frac{2(1-\alpha)}{z^2}+\sum_{k=1}^{\infty}(k-1+2\alpha)|\sigma_k(\alpha_1)|a_k z^{k-1}} \right| < \eta.$$

Since $|\text{Re}(z)| \le |z|$ for all z , we get

$$Re\left\{ \frac{\sum_{k=1}^{\infty}(k+1)|\sigma_k(\alpha_1)|a_k z^{k-1}}{\frac{2(1-\alpha)}{z^2}+\sum_{k=1}^{\infty}(k-1+2\alpha)|\sigma_k(\alpha_1)|a_k z^{k-1}} \right\} < \eta \quad (7)$$

on the real axis when choosing the value of $z$ thus the value of

$$\frac{z\,(W(\alpha_1)f(z)\,)''}{(W(\alpha_1)f(z))'}\ ,$$

is real , therefore clearing the denominator of (7) and when z →1⁻ through real axis the result is sharp for the function

$$f_k(z) = \frac{1}{z} + |\sigma_k\,(\alpha_1)|^{-1} \times \frac{2\eta(1-\alpha)}{[k(1+\eta)+(1+\eta(2\alpha-1))]} z^k, \text{k} \ge 1 \quad (8)$$

**Corollary 2.1** : When $f \in V(\alpha,\eta)$, then

$$a_k \le |\sigma_k(\alpha_1)|^{-1} \times \frac{2\eta(1-\alpha)}{[k(1+\eta)+(1+\eta(2\alpha-1))]},$$

where $0 < \alpha < 1$ , $0 < \eta \le 1$ .

## 3. Growth and Distortion Theorems

Distortion and growth Theorems property for the function $f \in V(\alpha,\eta)$, is given as follows :

**Theorem 3.1:** Let $f$ be a function defined by (1) is in the class $V(\alpha,\eta)$.
Then for $0 < |z| = r < 1$ we get

$$\frac{1}{r} - r|\sigma_1(\alpha_1)|^{-1} \times \frac{\eta(1-\alpha)}{(1+\alpha\eta)} \le |f(z)| \le \frac{1}{r} + r\frac{\eta(1-\alpha)}{(1+\alpha\eta)}|\sigma_1(\alpha_1)|^{-1},$$

equivalences for

$$f(z) = \frac{1}{z} + |\sigma_1(\alpha_1)|^{-1} \times \frac{\eta(1-\alpha)}{(1+\alpha\eta)} z.$$

**Proof:** Since $f \in V(\alpha,\eta)$, then we get by theorem 2.1, then the inequality

$$\sum_{k=1}^{\infty}\sigma_k(\alpha_1)[k(\eta+1)+(1+\eta(2\alpha-1))]a_k \le 2\eta(1-\alpha).$$

Then

$$|f(z)| \le \left|\frac{1}{z}\right| + \sum_{k=1}^{\infty} a_k |z|^k,$$

for $0 < |z| = r < 1$ , we get

$$|f(z)| < \frac{1}{r} + r\sum_{k=1}^{\infty} a_k \le \frac{1}{r} - |\sigma_1(\alpha_1)|^{-1} \times \frac{\eta(1-\alpha)}{(1+\alpha\eta)} r.$$

In addition to
$$|f(z)| \ge \left|\frac{1}{z}\right| - \sum_{k=1}^{\infty} a_k |z|^k \ge \frac{1}{r} - |\sigma_1(\alpha_1)|^{-1} \times \frac{\eta(1-\alpha)}{(1+\alpha\eta)} r, \ |z| = r.$$

**Theorem 3.2:** Let A function $f$ defined by (1) in the class $f \in V(\alpha,\eta)$. Then *for* $0 < |z| = r < 1$ we get

$$\frac{1}{r^2} - |\sigma_1(\alpha_1)|^{-1} \times \frac{\eta(1-\alpha)}{(1+\alpha\eta)} \leq |f'(z)| \leq \frac{1}{r^2} + |\sigma_1(\alpha_1)|^{-1} \times \frac{\eta(1-\alpha)}{(1+\alpha\eta)}.$$

Equivalences for

$$f(z) = \frac{1}{z} + |\sigma_1(\alpha_1)|^{-1} \times \frac{\eta(1-\alpha)}{(1+\alpha\eta)} z.$$

**proof:** Form Theorem 2.1, we get

$$\sum_{k=1}^{\infty} |\sigma_k(\alpha_1)|[k(1+\eta) + (1+\eta(2\alpha-1))]a_k \leq 2\eta(1-\alpha).$$

Thus

$$|f'(z)| \leq \left|\frac{-1}{z^2}\right| + \sum_{k=1}^{\infty} ka_k |z|^{k-1},$$

for $0 < r = |z| < 1$ we get

$$|f'(z)| \leq \left|\frac{-1}{r^2}\right| + \sum_{k=1}^{\infty} ka_k$$

$$\leq \frac{1}{r^2} + |\sigma_1(\alpha_1)|^{-1} \times \frac{\eta(1-\alpha)}{(1+\alpha\eta)}.$$

And

$$|f'(z)| \geq \left|\frac{-1}{z^2}\right| - \sum_{k=1}^{\infty} ka_k |z|^{k-1},$$

$$\geq \left|\frac{1}{r^2}\right| - \sum_{k=1}^{\infty} ka_k$$

$$\geq \frac{1}{r^2} - |\sigma_1(\alpha_1)|^{-1} \times \frac{\eta(1-\alpha)}{(1+\alpha\eta)}.$$

### 4. Hadamard product

**Theorem 4.1:** If the function $g, f \in V(\alpha, \eta)$. Then $(f * g) \in V(\alpha, \eta)$, for

$$f(z) = \frac{1}{z} + \sum_{k=1}^{\infty} a_k z^k,$$

$$g(z) = \frac{1}{z} + \sum_{k=1}^{\infty} b_k z^k,$$

and

$$(f * g)(z) = \frac{1}{z} + \sum_{k=1}^{\infty} a_k b_k z^k,$$

where

$$\delta = \frac{2\eta^2 (1-\alpha)(k+1)}{2\eta^2(1-\alpha)(k+2\alpha-1) - |\sigma_n(\alpha_1)|[k(1+\eta) + (1+\eta(2\alpha-1))]^2}$$

**Proof:** Since $f, g \in V(\alpha, \eta)$, then by Theorem 2.1, we have

$$\sum_{k=1}^{\infty} |\sigma_k(\alpha_1)| \frac{[k(1+\eta) + (1+\eta(2\alpha-1))]}{2\eta(1-\alpha)} a_k \leq 1,$$

and

$$\sum_{k=1}^{\infty} |\sigma_k(\alpha_1)| \frac{[n(1+\eta) + (1+\eta(2\alpha-1))]}{2\eta(1-\alpha)} b_k \leq 1,$$

*we must find the largest δ such that*

$$\sum_{k=1}^{\infty}|\sigma_k(\alpha_1)|\frac{[k(1+\eta)+(1+\eta(2\alpha-1))]}{2\eta(1-\alpha)}\sqrt{a_k b_k} \leq 1. \qquad (9)$$

To prove the theorem it is over to show that

$$|\sigma_k(\alpha_1)|\frac{[k(1+\delta)+(1+\delta(2\alpha-1))]}{2\delta(1-\alpha)}a_k b_k$$
$$\leq |\sigma_k(\alpha_1)|\frac{[k(1+\eta)+(1+\eta(2\alpha-1))]}{2\eta(1-\alpha)}\sqrt{a_k b_k} \quad ,$$

which is equivalent to

$$\sqrt{a_k b_k} \leq \frac{\delta[k(1+\eta)+(1+\eta(2\alpha-1))]}{\beta[k(1+\delta)+(1+\delta(2\alpha-1))]} \quad .$$

From (9) we get

$$\sqrt{a_k b_k} \leq |\sigma_r(\alpha_1)|\frac{2\eta(1-\alpha)}{[k(1+\eta)+(1+\eta(2\alpha-1))]} \quad .$$

We must proof that

$$|\sigma_k(\alpha_1)|\frac{2\eta(1-\alpha)}{[k(1+\eta)+(1+\eta(2\alpha-1))]} \leq \frac{\delta[k(1+\eta)+(1+\eta(2\alpha-1))]}{\eta[k(1+\delta)+(1+\delta(2\alpha-1))]},$$

which gives

$$\delta \leq \frac{2\eta^2(\alpha-1)(k+1)}{2\eta^2(1-\alpha)(k+2\alpha^{-1})-|\sigma_k(\alpha_1)|[k(1+\eta)+(1+\eta(2\alpha-1))]^2} \quad .$$

**Theorem 4.2** : If the function $f_i$ (i=1 ,2) defined by

$$f_i(z) = \frac{1}{z} + \sum_{k=1}^{\infty} a_{k,i} z^k \, , (a_{k,i} \geq 0 \, , i = 1,2)$$

be in the class $V(\alpha,\eta)$, then the function defined

$$g(z) = \frac{1}{z} + \sum_{k=1}^{\infty}(a^2_{k,1} + a^2_{k,2})z^k \, ,$$

is in the class $V(\alpha,\eta)$, where

$$\beta = \frac{4\eta^2(\alpha-1)(k+1)}{4\eta^2(\alpha-1)(k+2\alpha-1)-|\sigma_k(\alpha_1)|[k(1+\eta)+(1+\eta(2\alpha-1))]^2}$$

**proof:** Since $f_i \in V(\alpha,\eta)$, (i= 1 ,2), then by Theorem 2.1, we get

$$\sum_{k=1}^{\infty}|\sigma_k(\alpha_1)|\frac{[k(1+\eta)+(1+\eta(2\alpha-1))]}{2\eta(1-\alpha)}a_{k,i} \leq 1, (i=1,2) \, .$$

**Hence**

$$\sum_{k=1}^{\infty}(|\sigma_k(\alpha_1)|\frac{[k(1+\eta)+(1+\eta(2\alpha-1))]}{2\eta(1-\alpha)})^2 a_{k,i}^2$$
$$\leq (\sum_{k=1}^{\infty}|\sigma_k(\alpha_1)|\frac{[k(1+\eta)+(1+\eta(2\alpha-1))]}{2\eta(1-\alpha)}a_{k,i})^2 \leq 1, (i=1,2).$$

Thus

$$\sum_{k=1}^{\infty} \frac{1}{2}|\sigma_k(\alpha_1)|(\frac{[k(1+\eta)+(1+\eta(2\alpha-1))]}{2\eta(1-\alpha)})^2(a_{k,1}{}^2 + a_{k,2}{}^2) \leq 1,$$

to prove the theorem we must find the largest $\beta$ such that

$$\frac{[k(\beta+1)+(1+\beta(2\alpha-1))]}{\beta} \leq \frac{|\sigma_k(\alpha_1)|[k(1+\eta)+(1+\eta(2\alpha-1))]^2}{4\eta^2(1-\alpha)}, k \geq 1,$$

so that

$$\beta \leq \frac{4\eta^2(\alpha-1)(k+1)}{4\eta^2(\alpha-1)(k+2\alpha-1)-|\sigma_k(\alpha_1)|[k(1+\eta)+(1+\eta(2\alpha-1))]^2}.$$

**Theorem 4.3:** If $f(z) = \frac{1}{z} + \sum_{k=1}^{\infty} a_k z^k \in V(\alpha,\eta)$, and

$g(z) = \frac{1}{z} + \sum_{k=1}^{\infty} b_k z^k$ with $|b_k| \leq 1$ is in the class $V(\alpha,\eta)$
then $f(z)*g(z) \in V(\alpha,\eta)$.

**Proof**: By Theorem 2.1, we get

$$\sum_{k=1}^{\infty} |\sigma_k(\alpha_1)|[k(1+\eta)+(1+\eta(2\alpha-1))]a_k \leq 2\eta(1-\alpha).$$

Since

$$\sum_{k=1}^{\infty} |\sigma_k(\alpha_1)|(\frac{[k(1+\eta)+(1+\eta(2\alpha-1))]}{2\eta(1-\alpha)}|a_k b_k|,$$

$$= \sum_{k=1}^{\infty} |\sigma_k(\alpha_1)|(\frac{[k(1+\eta)+(1+\eta(2\alpha-1))]}{2\eta(1-\alpha)} a_k |b_k|,$$

$$\leq \sum_{k=1}^{\infty} |\sigma_k(\alpha_1)|[k(1+\eta)+(1+\eta(2\alpha-1))]a_k \leq 1.$$

Thus $f(z)*g(z) \in V(\alpha,\eta)$.

**Corollary 4.1:** If $f(z) = \frac{1}{z} + \sum_{k=1}^{\infty} a_k z^k \in V(\alpha,\eta)$, and $g(z) = \frac{1}{z} + \sum_{k=1}^{\infty} b_k z^k$ with $0 \leq b_k \leq 1$ is in the $V(\alpha,\eta)$, then $f(z)*g(z) \in V(\alpha,\eta)$.

## 5. Radil of starlikness and convexity

**Theorem 5.1**: Let $f(z)$ be the function defined by (1) be in the subclass $V(\alpha,\eta)$. Then f is meromorphically starlike of order $\delta$ ($0 \leq \delta < 1$) in the disk $|z| < r_1(\alpha,\eta,\delta)$, where

$$r_1(\alpha,\eta,\delta) = \inf_k \{|\sigma_k(\alpha_1)| \frac{[k(1+\eta)+(1+\eta(2\alpha-1))](1-\delta)}{2\eta(k+2-\delta)(1-\delta\alpha)}\}^{\frac{1}{k+1}}$$

the result is sharp for the function given by (8).

**Proof:** We show that

$$\left|\frac{zf'(z)}{f(z)}+1\right| \leq 1-\delta,$$

$$\left|\frac{zf'(z)}{f(z)}+1\right| = \left|\frac{\sum_{k=1}^{\infty}(k+1)a_k z^k}{z^{-1}+\sum_{k=1}^{\infty}a_k z^k}\right| \leq \frac{\sum_{k=1}^{\infty}(k+1)a_k |z|^{k+1}}{1-\sum_{k=1}^{\infty}a_k |z|^{k+1}}.$$

This will be bounded by $1-\delta$,

$$\frac{\sum_{k=1}^{\infty}(k+1)a_k |z|^{k+1}}{1-\sum_{k=1}^{\infty}a_k |z|^{k+1}} \leq 1-\delta,$$

$$\sum_{k=1}^{\infty}(2+k-\delta)a_k |z|^{k+1} \leq 1-\delta,$$

from Theorem 2.1, we get

$$\sum_{k=1}^{\infty}|\sigma_k(\alpha_1)|\left(\frac{[k(1+\eta)+(1+\eta(2\alpha-1))]}{2\eta(1-\alpha)}\right)a_k \leq 1.$$

Hence

$$|z|^{k+1} \leq |\sigma_k(\alpha_1)|\frac{[k(1+\eta)+(1+\eta(2\alpha-1))](1-\delta)}{2\eta(k+2-\delta)(1-\alpha)},$$

$$|z| \leq \left\{|\sigma_k(\alpha_1)|\frac{[k(1+\eta)+(1+\eta\ (2\alpha-1))](1-\delta)}{2\eta(k+2-\delta)(1-\alpha)}\right\}^{\frac{1}{k+1}} \ \frac{1}{k+1}.$$

**Theorem 5.2:** Let the function $f(z)$ defined by (1) be in the subclass $V(\alpha,\eta)$. Then $f$ is meromorphically convex of order $\gamma$ $(0\leq\gamma<1)$ in the disk $|z|<r_2(\eta,\alpha,\gamma)$, where

$$r_2(\eta,\alpha,\gamma)= \inf\left\{|\sigma_k(\alpha_1)|\frac{[k(1+\eta)+(1+\eta(2\alpha-1))](1-\gamma)}{2\eta(k+2-\gamma)(1-\alpha)}\right\}^{\frac{1}{k+1}}.$$

The result is sharp for the function given by (7).

**Proof:** By utilizing the same way in the proof of theorem 5.1 we can get this

$$\left|\frac{zf''(z)}{f'(z)}+2\right| \leq 1-\gamma, ((0\leq\gamma<1)).$$

For $|z|< r_2$ depending on the help of the Theorem 2.1, lead to confirmed of theorem 5.2.□